# Generalization of Cantor Pairing Polynomials (Bijective Mapping Among Natural Numbers) from $N_0^2 \to N_0$ to $Z^2 \to N_0$ and $N_0^3 \to N_0$

Sandor Kristyan

*Institute of Materials and Environmental Chemistry, Research Centre for Natural Sciences,
Eötvös Loránd Research Network (ELKH), H-1117 Budapest, Magyar Tudósok körútja 2, Hungary*

Corresponding author: kristyan.sandor@ttk.mta.hu, kristyan.sandor@ttk.hu

**Abstract**. The Cantor pairing polynomials are extended to larger 2D sub-domains and more complex mapping, of which the most important property is the bijectivity. If "corners" are involved inside (but not the borders of) domain, more than one connected polynomials are necessary. More complex patterns need more complex subsequent application of math series to obtain the mapping polynomials which are more and more inconvenient, although elementary. A tricky polynomial fit is introduced (six coefficients are involved like in the original Cantor polynomials with rigorous but simple restrictions on points chosen) to buy out the regular treatment of math series to find the pairing polynomials instantly. The original bijective Cantor polynomial $C_1(x,y) = (x^2+2xy+y^2+3x+y)/2$: $N_0^2 \to N_0$ ($\equiv$ positive integers) which is 2-fold and runs in zig-zag way along lines $x+y=N$ is extended e.g. to the bijective $P(x,y) = 2x^2+4\text{sgn}(x)\text{sgn}(y)xy+2y^2-2H(x)\text{sgn}(y)x-y+1$: $Z^2 \to N_0$ (with sign and Heaviside functions, $Z \equiv \pm$integers) running in spiral way along concentric rhombuses, or to the bijective $P_{3D}(x,y,z) = [x^3+y^3+z^3 +3(xz^2+yz^2+zx^2+2xyz+zy^2+yx^2+xy^2) +3(2x^2+2y^2+z^2+2xz+2yz+4xy) +5x+11y+2z]/6$: $N_0^3 \to N_0$ which is 6-fold and runs along plains $x+y+z=N$. Storage device for triangle matrices is also commented as "cutting" the original Cantor domain to half along with related Diophantine equations.

**Keywords**. Generalization of Cantor pairing polynomials; Extension of Cantor's 2D domain to $Z^2 \to N_0$ and $N_0^3 \to N_0$

## INTRODUCTION

A pairing function is a process to uniquely encode two natural numbers into a single natural number. It is a bijection $N_0^2 \equiv N_0 \times N_0 \to N_0$, where $N_0 \equiv \{0,1,2,3,\ldots\}$ is the set of natural numbers with zero i.e. the nonnegative integers. $Z = \{0, \pm1, \pm2, \pm3,\ldots\}$ notates the integer numbers. It is also called packing polynomial that maps the lattice points of $N_0^2$ with nonnegative coordinates bijectively onto $N_0$. Cantor constructed two (but essentially the same) quadratic packing polynomials. The so-called Cantor polynomial(s) [1]

$$C_1(x,y) = \tfrac{1}{2}(x+y)^2 + \tfrac{1}{2}(3x+y) = \tfrac{1}{2}(x^2+2xy+y^2+3x+y) = \binom{x}{1} + \binom{x+y+1}{2} \quad (1)$$

$$C_2(x,y) = C_1(y,x) = \tfrac{1}{2}(x+y)^2 + \tfrac{1}{2}(3y+x) = \tfrac{1}{2}(x^2+2xy+y^2+3y+x) = \binom{y}{1} + \binom{x+y+1}{2} \quad (2)$$

is a bijective mapping $N_0^2 \to N_0$; the far right is the more compact combinatorial form. Fueter and Pólya stated the theorem [2] that, there are no other quadratic polynomials with this property. The geometric picture on Figure 1 is quite impressive how this mapping works. In reverse way, if $C_1(x,y) = z \in N_0$ is given, then with $w \equiv \text{int}(-\tfrac{1}{2} + \tfrac{1}{2}\text{sqrt}(1+8z))$ the inverse of Eq.1 is

$$C_1^{-1}(z) = (x := z-w(w+1)/2,\ y := w(w+3)/2-z) = (x := z-w(w+1)/2,\ w-x) . \quad (3)$$

Obviously, the two polynomials in Eqs.1-2 are mirror image of each other with respect to line $y=x$ in the first quarter (+,+) in $Z^2$, i.e. in $N_0^2$. The proof of Eqs.1-2 is elementary [3], in fact it is based on the subsequent use of sum of mathematical series ($\Sigma_{i=1}^{n} a_i = (2a_1+(n-1)d)n/2$, e.g. $\Sigma_{i=1}^{x} i = (1+x)x/2$ with $x \in N_0$), the latter yields that Eqs.1-2 are second order polynomials.

In addition to Cantor polynomials, a few other pairing functions are also described in the literature. Most notably, the Rosenberg-Strong pairing function [4] for the non-negative integers is defined, see below. The practical use of these pairing functions is benefited for example, in computer sciences, where there is the problem of saving and retrieving matrices (their indices $\in N_0^2$) and other multidimensional data in linear memory (their indices $\in N_0$). Theoretically, it is also interesting to consider the Cantor polynomials more generally, for example, the cases indicated in the title. Notations are listed in App.1.

## SIMPLE MANIPULATIONS

Next, we will deal with different patterns and domains than the $N_0^2$ on Fig.1, where the mirror image property is also present, but we will consider only one of these mirror image(s) for the sake of brevity. The bijective property allows to recompose the Cantor polynomials in relation to Diophantine equations:

**Lemma 1.** For $a,b,z \in N_0$ and z is given, Eq.3 provides the only one solution for Diophantine equation

$$2z = a^2 + 2ab + b^2 + 3a + b = (a+b)^2 + 3a + b . \quad (4)$$

If z is not given, then for any two solutions $(a_i,b_i,z_i)$ of Eq.4 with i=1 and 2 in which $(a_1,b_1) \neq (a_2,b_2)$, the $z_1 \neq z_2$.

The first part can be pictured or clear via Fig.1 (both cases). For example, if z=7 is given, then (a,b)=(1,2), as well as any other value, e.g. a=2 in Eq.4 yields $b^2+5b-4=0$ yielding b= ±½sqrt(41)-5/2 non-integer solution, generally $b^2+(2a+1)b+a^2+3a-2z=0$ with b= ±½sqrt(1+8(z-a)) –a–½. If z is not given, then any (a,b) for Eq.4 gives a defined even number, that is, if $(a,b) \neq (a',b')$ then $2z \neq 2z'$. The following two lemmas are elementary and hold for all extensions of Cantor polynomials especially for the discussed cases in this work:

**Lemma 2.** The linear transformation of the image set and shift of the domain are as follows: The transformation $k_1C_i(x,y)+k_2$ in Eqs.1-2 maps $N_0^2$ to 1D array $\{k_2, k_1+k_2, 2k_1+k_2, 3k_1+k_2,…\}$, where $k_1,k_2 \in Z$, but can be real or complex as well.

For example with $(k_1,k_2)=(2,1)$ the original Cantor image set $N_0 \equiv \{0,1,2,3,4,…\}$ transforms to positive odd numbers $\{1,3,5,7,9,…\}$. The transformation $C_i(x-k_1,y-k_2)$, where $k_1,k_2 \in Z$ keeps the image set $N_0$, but shifts the $N_0^2$ domain (Fig.1) with vector $(k_1,k_2)$, for example $(k_1,k_2)=(0,1)$ on (both in) Fig.1 shifts the lattice points upward with 1 (e.g. on left of Fig.1 the y coordinates of numbers 0,2,5 and 9 change from 0 to 1); notice that in the shift of domain (except certain subsets at boundary) the new lattice points overlap with the old ones.

**Lemma 3.** Rotations (around perpendicular z-axis to (x,y)-plain on Fig.1) are possible. For example, Eq.1 plotted on left of Fig.1 can be rotated to clockwise with π/2 (e.g. image value 2 changes position as (x=1,y=0)→(X=0, Y=-1)), plotted on left of Figure 2. (Notice that rotations with ±π/2, ±π and ±3π/2 around the origin make complete overlap between the grid and its image, otherwise not.) Algebraically, this rotation is a standard linear transformation and re-notation as (x,y) = (-Y,X) → (-y,x), and, for example, Eq.1 transform as

$$C_{1rot}(x,y) = \tfrac{1}{2}(x^2 - 2xy + y^2 + x – 3y) . \quad (5)$$

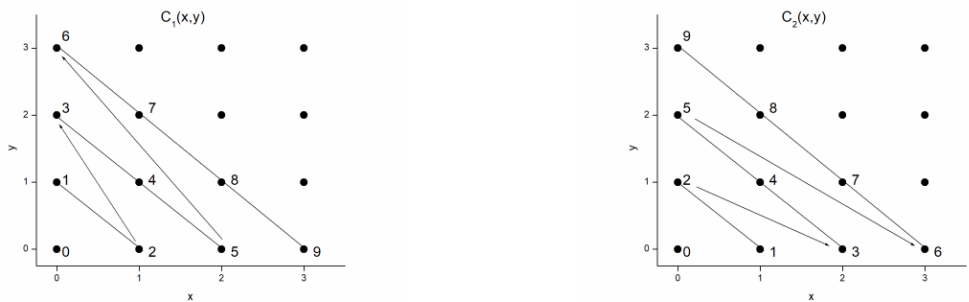

**FIGURE 1.** Bijective mapping of Cantor polynomials: $N_0^2 \to N_0$ in Eqs.1-2, e.g. $C_1(1,0)=2$, $C_2(1,0)=1$, $C_1(1,1)=C_2(1,1)=4$, etc.. The mapping $C_1$ [Left] and $C_2$ [Right] are mirror images to each other with respect to line y=x.

One can cut out the hectic treatment of math series with a simple linear algebraic fitting procedure to find the Cantor-like polynomials. Eqs.1-2 are 2nd order polynomials of the form (we avoid double indices for simplicity)

$$P(x,y) = a_6 x^2 + a_5 xy + a_4 y^2 + a_3 x + a_2 y + a_1 \quad (6)$$

with zero constant term ($a_1=0$, although $a_1 \neq 0$ can also come up e.g. by Lemma 2). Using Eq.6, the Eq.5 can be derived in another way than above. Knowing that the pairing function is a unique second order polynomial (Fueter-Pólya) and the rotation preserves this form, fitting (more exactly fixing) polynomial in Eq.6 is possible with using an adequate point set with cardinality six, e.g. set A on left of Fig.2, and this P(x,y) will work for all other lattice points in the domain. This adequate point set (of six) is the "initial triangle part" of the domain, like points marked with letter A on left of Fig.2, i.e. with function values {0,1,2,3,4,5}. Another (i.e. any non-linear which avoids zero determinant) point set is not adequate generally, its pattern must "represent" the domain, regardless that, generally, a polynomial fit (here fix) for its point set (containing as much as many parameters it possesses, e.g. 6 in Eq.6) is unique. The reason for this peculiarity is detailed below and exampled with Eq.9 with right of Figure 3 and Eq.10 with right of Fig.2 below.

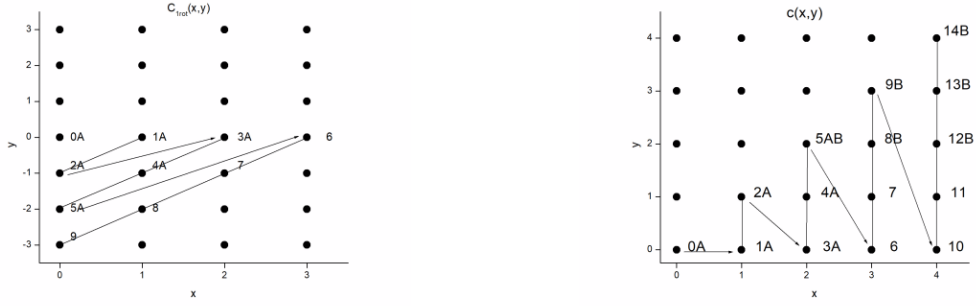

**FIGURE 2.** Clockwise rotation by $\pi/2$ of mapping $C_1(x,y)$ on [Left]: $C_{1rot}$: $\{(x,y)| x\geq 0, y\leq 0, x\in N_0, -y\in N_0\}\to N_0$ in Eq.5. Letter A marks the six points (the "initial triangle part"), for example, which can be used by algorithm in Eq.7 to obtain Eq.5 [Left]. Bijective mapping of $c(x,y)$: $\{(x,y)| 0\leq y\leq x, x\in N_0, y\in N_0\}\subset N_0^2\to N_0$ in Eq.10 (i.e. between y=0 and y=x in the (+,+) quarter). Any (non-linear) set of six points cannot guarantee the right operation of algorithm in Eq.7, but set A= $\{0,1,2,3,4,5\}$ (the "initial triangle part") in the neighborhood of origin and set B= $\{5,8,9,12,13,14\}$ marked with letters A and B resp. provide the right pairing polynomial by algorithm in Eq.7 resulted in Eq.10 [Right].

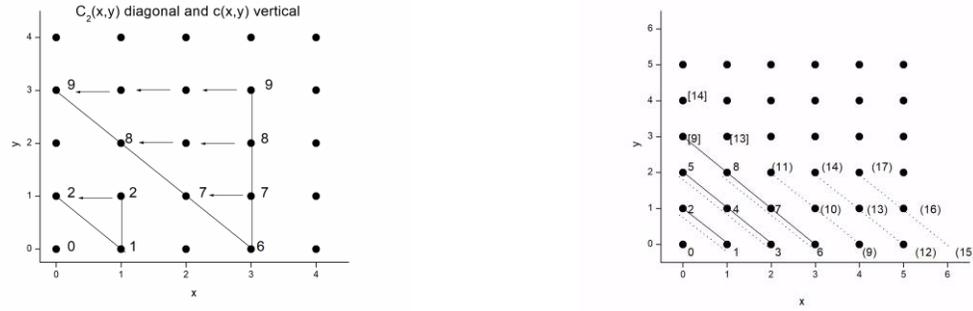

**FIGURE 3.** Horizontal shift/transformation between $C_2(x,y)$ diagonal (Eq.2 and right of Fig.1) as well as $c(x,y)$ vertical (Eq.10 and right of Fig.2) values [Left]. Cantor mapping (solid lines right of Fig.1) and saw mapping (dotted lines) on lattice points. Common points values are without brackets, the subset in both as $\{0,1,2,3,4,5,6,7,8\}$, Cantor point values not overlapping with saw points are in square brackets, saw point values not overlapping with Cantor points are in brackets [Right].

Generally, matrix $[A_{6\times 6}]$ with rows $[x_i^2, x_i y_i, y_i^2, x_i, y_i, 1]$ for i=1,2,3,4,5,6 must be set up for the adequate point set with function values $P(x_i,y_i)$ and one has to solve the linear equation system

$$[A_{6\times 6}] [a_6,a_5,a_4,a_3,a_2,a_1]^T = [P(x_1,y_1), P(x_2,y_2), P(x_3,y_3), P(x_4,y_4), P(x_5,y_5), P(x_6,y_6)]^T . \quad (7)$$

We will use this simple algorithm in Eq.7 for all the more complex domains below to obtain the parameters for Cantor-like (or extended Cantor) polynomials or pairing equations. Particularly, algorithm in Eq.7 yields the simple linear equation system for the coefficients in Eq.5 if one substitutes $(x,y,C_{1rot}(x,y))=(0,0,0), (1,0,1), (0,-1,2), (2,0,3), (1,-1,4)$ and $(0,-2,5)$ marked with letter A on left of Fig.2 into Eqs.6-7. As expected, this equation system yields the solution $[½, -1, ½, ½, -3/2, 0]^T$ which is the same as was obtained in Eq.5 with rotation transformation discussed; notice that this point set is the "initial triangle part" of the domain mentioned above. For example, point set of six with function values $\{1,2,6,7,8,9\}$ on left of Fig.2 yields the wrong polynomial

$$P_{wrong}(x,y)= 0.75 \, x^2 -1.5 \, xy +0.75 \, y^2 -0.5x -0.5 \, y +0.75 \quad (8)$$

which reproduces the six points used in the fit, but other points are not reproduced, e.g. P(0,0)= 0.75 instead of 0 or P(4,0)= 10.75 instead of 10. Of course, Eqs.1-2 for Fig.1 can also be derived with the procedure in Eqs.6-7. It will be seen below that, the algorithm or procedure in Eq.7, e.g. for Figs.1-2a to obtain Eqs.1,2 and 5, is much simpler and transparent, than for example re-doing the (original) derivation (as used for left of Fig.1 to obtain Eq.1 in the literature) with applying subsequent mathematical series. If a $P(x,y)$ is obtained by the algorithm in Eq.7 for a mapping (like the ones above for Figs.1-2a with using a chosen point set of six), it can be tested immediately on the figure if it works (or not coming from the wrong point set chosen for the fit).

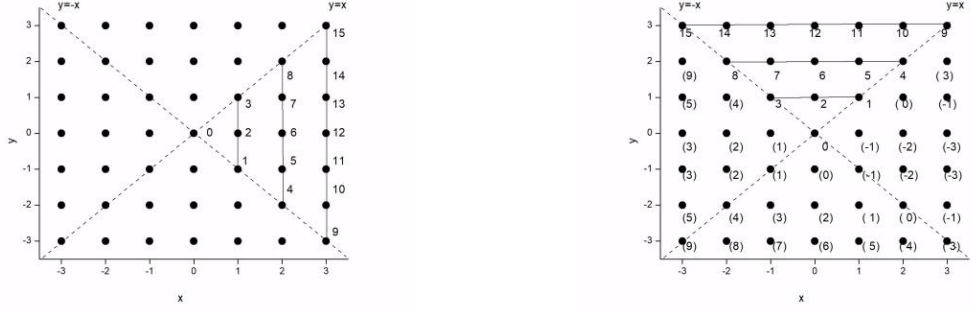

**FIGURE 4.** Bijective mapping (triangle) of $P_{\text{left of Fig.4}}(x,y)$: $\{(x,y)|\ -x\leq y\leq x,\ x\in N_0,\ y\in Z\}\to N_0$ in Eq.13 (i.e. between $y=\pm x$ in the $(+,\pm)$ quarters) [Left]. Bijective mapping (triangle) of $P_{\text{right of Fig.4}}(x,y)$: $\{(x,y)|\ -y\leq x\leq y,\ x\in Z,\ y\in N_0\}\to N_0$ in Eq.13 (i.e. between $y=\pm x$ in the $(\pm,+)$ quarters). Values without brackets are the $P_{\text{right of Fig.4}}(x,y)$ with lattice points $(x,y)\in$ Cantor-domain of $P_{\text{right of Fig.4}}(x,y)$ defined here, the values in brackets are the $P_{\text{right of Fig.4}}(x,y)$ with lattice points $(x,y)$ in out-of-domain. Notice that the image or $P_{\text{right of Fig.4}}(x,y)$ values (in and not in brackets) are the 1D integer Z as horizontal "tiles", algebraically the $-x$ in Eq.13, the shift between adjacent tiles is $y(y+1)$, especially at $y=0$ and $-1$ the shift is $0(0+1) = -1(-1+1)= 0$ [Right].

**Definition 1.** We call "Cantor-like polynomials" (or extended Cantor polynomials) the extended pairing equations for larger and/or more complex domains than the domain (Fig.1) of original Cantor polynomials (Eqs.1-2). We call "algorithm (or procedure to solve the linear system) in Eq.7" to obtain the parameters for Cantor-like polynomials, this algorithm also includes a wise choice (see "initial triangle part" above and details next) of six points to fit.

There is one important restriction on algorithm in Eq.7. It is known that a polynomial fit is unique (e.g. two points determine a line, etc.), here the number of points to fit is the same as the number of parameters. Eq.6 is determined by six, linearly independent (that is, the determinant in Eq.7 is non-zero) points. However, now the domains brings special feature into the fit, e.g. on Fig.1 the points on x and y axis possess "breaks" in continuation to the subsequent function values ($n\to n+1$). It puts an extra requirement on the polynomial which is, in fact, independent of the fit which is unique itself. The consequence of these "breaks in the mapping" can be displayed simply with the example on right of Fig.3, where two mappings, the Cantor (right of Fig.1) and a "saw mapping" are plotted together. One can see that the location of function values $\{0,1,2,3,4,5,6,7,8\}$ are common in both mappings, but other $(x,y)$ coordinates yield different function values, e.g. $(x=2,y=2)$ Cantor mapping/polynomial yields [12] (not marked), but saw mapping/polynomial yields (11) marked. Obviously a set containing 6 points from that common set yields one certain polynomial for Eq.6 by the algorithm in Eq.7 which cannot work for both mappings at the same time on right of Fig.3. It means that wrong choice of 6 points (to fit the 6 parameters in Eq.6) yields wrong polynomial for Eq.6 when one wants to describe any of the two mappings on right of Fig.3. For example, function values $\{0,1,2,3,4,5\}$ with its $(x,y)$ coordinates for Eq.6 yields Eq.2 by the algorithm in Eq.7, this set for fit correctly represents the periodicity, however it does not yield correct polynomial for saw mapping, because this set does not represent the periodicity of that. Set with function values $\{0,1,2,3,4,6\}$ with its $(x,y)$ coordinates for Eq.6 yields zero determinant by the algorithm in Eq.7 and does not provide any polynomial. Also, point set with function values $\{1,2,3,4,6,7\}$ yields $P(x,y)=\frac{1}{2}x^2+xy+4.6915y^2+\frac{1}{2}x-2.6915y$, a wrong polynomial (like Eq.8) which fails for both mappings on right of Fig.3, because this point set does not represent the periodicity of any of the two mappings plotted. However, set with function values $\{3,6,7,(9),(10),(11)\}$ with its $(x,y)$ coordinates for Eq.6 yields $P(x,y)=3x+4y-3$ by the algorithm in Eq.7 for saw mapping correctly, because the set represents the periodicity correctly (also an "initial triangle" part). One must notice that, for saw mapping, there are two domains on right of Fig.3, domain I contains $(x,y)$ lattice points with $(0\leq x\leq 1, 0\leq y\leq 2)$, here Eq.2 holds, the Cantor $2^{nd}$ order polynomial, as well as domain II contains $(x,y)$ lattice points with $(x>1, 0\leq y\leq 2)$, here the above $1^{st}$ order polynomial holds obtained above, altogether

$$P_{\text{right of Fig.3}}(0\leq x\leq 1, 0\leq y\leq 2) = C_2(x,y) \equiv \tfrac{1}{2}(x^2 + 2xy + y^2 + 3y + x) \quad \text{and} \quad P_{\text{right of Fig.3}}(1<x, 0\leq y\leq 2) = 3x+4y-3. \quad (9)$$

In summary from the examples above, the pattern/shape/topology of the point set of six for Eq.7 is fundamental. Notice that saw mapping needs only first order polynomial on its domain II. (By the way, the simplest mapping is the simplest saw mapping, that is the x axis to the x axis itself where $(x,y=0)\in N_0^2 \to x\in N_0$, a $1^{st}$ order polynomial, $P(x,y=0)=x$.) Another derivation (than the algorithm in Eq.7) for domain II for saw mapping is as follows. The $P_{\text{right of Fig.3}}(1<x,y=0)= 3(x-1)$ what can be checked easily, and the vertical increases is +4 everywhere (e.g. at $x=3$ the 6+4 =(10) and (10)+4=(14)), so $P_{\text{right of Fig.3}}(1<x,y)= 3(x-1)+4y$ as in Eq.9. The vertical width of saw domain (a strip) on right of Fig.3 is n=3 lattice point rows (0≤y≤2). The width with value n>3 is analogous (0≤y≤n-1) and $P_{\text{right of Fig.3}}$(0≤x≤n-2,y)= $C_2(x,y)$ a $2^{nd}$ order polynomial and $P_{\text{right of Fig.3}}$(n-2<x,y) is a $1^{st}$ order polynomial; the vertical borderline (x=n-3) between domain I and II shifts by 1 as n increases by 1, so as width n→∞ the domain I increases to +x direction (becoming the (+,+) quarter of Cartesian lattice points) and domain II disappears by shifting into infinity. Table 1 and App.1 list some more mappings of semi-infinite domains with linear relationships in relation to Cantor-like polynomials.

As a rule of thumb, to choose the six points for algorithm in Eq.7 to set up Cantor-like polynomials, to avoid the although elementary but hectic treatment of mathematical series, is as follows. In many cases quarters or eighth parts of Cartesian plain constitute the regions, compare domains on Figs.2a-2b to domains on Figs.2e-5, so the simplest choice for six points is the ones marked with letter A on Figs.2a and 2b. The following note is interesting: In this case the determinant value (noted as det[$A_{6x6}$]) is independent of the point set chosen near to origin, its value depends on its pattern. For example, for the two point sets marked with letter A and letter B on right of Fig.2, and the four regions on left of Figure 6 (the four point sets listed below), all yield det[$A_{6x6}$]=-4. The different pattern of point set (of six), e.g. on Figs.2e and 2f (the two point sets listed below), both yield det[$A_{6x6}$]=-8. Finally, the pattern determines the det[$A_{6x6}$] value not its location, it can be rotated or shifted. This invariant property is not detailed more here. Furthermore, even more important (and interesting) that pattern (with angle $45^0$) of set of six points on right of Fig.2 (marked with letter A or B) can be used for domain (with angle $90^0$) on left of Fig.2 in algorithm in Eq.7, see example below for left of Figure 7 region II with Eq.27, this comes from the "tile property" of Cantor-like polynomials (Theorem 1).

**Definition 2.** For the domain, image, and bijective mapping, instead of the rigorous algebraic way (with relations ⊂, ∈, ≤, etc.), the domain, image, and bijective mapping will be defined "simply" by the graphical (or visual) presentation on the figure discussed. For example, the domain ($N_0^2$) and the procedure of bijective mapping ($N_0^2 \to N_0$) are quite obvious for Cantor polynomials by looking at Fig.1.

A very simple case of bijective mapping c(x,y): {(x,y)| 0≤y≤x, x∈$N_0$, y∈$N_0$} ⊂ $Z^2 \to N_0$ is defined on right of Fig.2. It can be related, for example, to 1D arrangement of elements in triangular matrices, or the counting (the image) and triangle domain on right of Fig.2 comes up in e.g. prime number generation and counting k= 2xy+x+y [5-6]. The derivation for its pairing polynomial with mathematical series is as follows. The vertical arrangement allows simple transparent treatment with math series. P≡ c(x,y=0) values (on x axis) are c(1,0)= 1, c(2,0)= 3= 1+2, c(3,0)= 6= 1+2+3, c(4,0)= 10= 1+2+3+4, c(5,0)= 15= 1+2+3+4+5, etc.. This is a sum of mathematical series as c(x,0)= $S_x$=1+2+3+…+x= (1+x)x/2. Here comes the simple part, the increment on vertical y is the y coordinate itself as c(x,y)= x(x+1)/2 + y. These vertical values are added as the values of 0≤y≤x, so

**Lemma 4.** Shrinking the domain is possible. The bijective mapping c(x,y): {(x,y)| 0≤y≤x, x∈$N_0$, y∈$N_0$} ⊂ $N_0^2$ → $N_0$ defined on right of Fig.2 is

$$c(x,y) = x(x+1)/2 + y = \tfrac{1}{2}x^2 + \tfrac{1}{2}x + y \qquad (10)$$

The cardinality of domain is still countable-infinite like in case of original Cantor polynomials.

As a proof, it is still simple by manipulating with math series (as was described just before Lemma 4, a standard way as for Eqs.1-2 in the literature), but generally manipulation with math series becomes very hectic if the domain is more and more complex as the ones below. The algorithm in Eq.7 for Eq.6 is simple, for example, the set A on right of Fig.2 provides (x,y,c(x,y))= (0,0,0), (1,0,1), (1,1,2), (2,0,3), (2,1,4) and (2,2,5) yielding matrix via Eq.7 with rows (0,0,0,0,0,1), (1,0,0,1,0,1), (1,1,1,1,1,1), (4,0,0,2,0,1), (4,2,1,2,1,1), (4,4,4,2,2,1) and constant vector [0,1,2,3,4,5]$^T$. Its solution is [½,0,0,½,1,0]$^T$, the parameter vector indicated in Eq.6-7 and 10 as expected.

Notice that, the simplicity of procedure in Eq.7 is practically the same for all the more complex domains (see below). Case of Eq.8 above exampled a wrong point set of six providing wrong pairing polynomial for that case, the choice of adequate point set is fundamental. The chosen six points from domain must represent the entire (semi-)infinite 2D domain. Here we show another good point set of six which is not the "initial triangle part" of the domain, but also represent the domain adequately: On right of Fig.2 the point set of six with function value {5,8,9,12,13,14} marked with letter B and the above used point set of six with function value {0,1,2,3,4,5} marked with letter A (the "initial triangle part") both yields Eq.10 with the algorithm in Eq.7.

Continuing the calculations for Cantor-like polynomials, Lemma 2 can be used. For example, to shift the domain on right of Fig.2 diagonally up by unity and to transform image $N_0$ to positive odd numbers {1,3,5,7,9,…} via Eq.10 is 2c(x-1,y-1)+1= $x^2$-x+2y-1, in which (1,1)→1, (2,1)→3, (2,2)→5, etc., it buys out re-doing the algorithm in Eq.7.

**Corollary of Lemmas 2-4.** For existing polynomials, the (connected) polynomials for bijective mappings for more complex domains (e.g., Figs.3-5 below, where fault lines and corners are present) can be obtained via the simple algorithm in Eq.7.

**Lemma 5.** Elementary functional transformations are possible. After a (simpler) Cantor-like polynomial is derived (e.g. Eqs.1-2 or Eq.10) functional transformations can be applied as $f(c(x,y))$ instead of repeating the derivation (via math series or polynomial fit by algorithm in Eq.7) from the beginning. For example, $C_{1or2}(x,y)$-4 shifts the starting zero value from origin (0,0) to point (1,1), see Fig.1 as well as Lemma 2, or changing the image values $\{0,1,2,3,4\ldots,n,\ldots\}$ on right of Fig.2 to squares $\{0,1,4,9,16,\ldots,n^2,\ldots\}$ simply needs from Eq.10

$$(c(x,y))^2 = (\tfrac{1}{2}x^2 + \tfrac{1}{2}x + y)^2 = x^4/4 + x^2/4 + y^2 + \tfrac{1}{2}x^3 + x^2 y + xy, \quad (11)$$

showing that this mapping needs $4^{th}$ order polynomial (i.e. the degree of Eq.6 is inadequately low).

Notice an important and interesting internal symmetry between Figs.1b versus 2b in relation how the points are connected. Topologically, the diagonals in the former contain the same image numbers as the verticals of the latter, although the domain of the former is $N_0^2$ while the domain of the latter is $\{(x,y)|\ 0 \leq y \leq x,\ x \in N_0,\ y \in N_0\}$. Even more, we can state the following:

**Lemma 6.** The internal role of Cantor polynomials for $N_0^2 \rightarrow N_0$ in linear transformations for $N_0^2 \rightarrow N_0^2$) is as follows: Analogously to Cantor bijective mapping/pairing polynomial(s) for $N_0^2 \rightarrow N_0$ (i.e. Eqs.1 or 2, which are mappings between 2D and 1D), bijective mapping between 2D (left of Fig.1) and a "smaller" domain of 2D (e.g. right of Fig.2) can be defined, for example with Eqs.1 and 10 as

$$B: N_0^2 \rightarrow C_1(x,y) := c(x,y) \in N_0 \rightarrow \{(x,y)|\ 0 \leq y \leq x,\ x \in N_0,\ y \in N_0\} \text{ by } c^{-1}(x,y) \quad \text{wherein} \quad B(x,y)=(x+y,x), \quad (12)$$

i.e. the net transformation by B is first order (linear) polynomial in its both variables, in matrix form $b \equiv [(1,1)^T,(1,0)^T]$ with $b \cdot (x,y)^T = (x+y,x)^T$ where b is the matrix representation of linear transformation B.

For example, $C_1(1,1) = 4 = c(2,1)$ between left of Fig.1 and 2b, that is $B(1,1)=(2,1)$. Eqs.1-2 and 12 reveal that the second order polynomial in Eqs.1-2 is necessary for mapping between 2D and 1D, but linear relationship is enough for the mapping in Eq.12 between 2D and 2D, the latter is well known in linear algebra. Eq.3 looks quadratic (second order) polynomial at the first glimpse, but w is not polynomial of z, and similarly for $c^{-1}$, however subsequent use of $C_1$ and $c^{-1}$ produces a linear (first order) polynomial in Eq.12. (Alternatively, supposing form $B:(x,y) \rightarrow (a_1 x + b_1 y, a_2 x + b_2 y)$, use the mapping of two points (with non-parallel position vector), for example, $(0,1) \rightarrow (1,0)$ and $(1,0) \rightarrow (1,1)$ and compare to the action of B as $(1,0)=(b_1,b_2)$ and $(1,1)=(a_1,a_2)$ providing the mapping in Eq.12.) Another representative example for Lemma 6 is between Figs.1a and 1b. The analogous mapping to Eq.12 is $B': N_0^2 \rightarrow N_0^2$, and the $C_1(x,y) := C_2(y,x)$ in Eq.2 yields $B'(x,y)=(y,x)$, that is, linear relation in both coordinates what can easily be verified on Fig.1, as well as the matrix representation of B' is $b'=[(0,1)^T,(1,0)^T]$.

The vertical, horizontal, and diagonal transformations are as follows: Interesting transformations can be applied to make the derivations of pairing polynomials simple. If increasing/decreasing function/image values are lying on vertical lines, the derivation with math series is simpler. Compare the literature derivation for Eqs.1-2 (Fig.1) to derivation here for Eq.10 (right of Fig.2), the vertical structure makes the derivation simpler for the latter as well as the derivation with algorithm in Eq.7 is the simplest. However, by horizontal shift, Eq.10 can be transformed to Eq.2. left of Fig.3 shows how right of Fig.1 and right of Fig.2 relate to each other. That is a vertical stretching by $\sqrt{2}$ (as $(x,y) \rightarrow (x,y\sqrt{2})$) and rotation of $(x,y)$ around $(x,0)$ by $-45°$ (sin or $\cos(45°)= 1/\sqrt{2}$ in rotation matrix annihilates the $\sqrt{2}$ from stretching to put beck the values to integers), or simply $(x,y) \rightarrow (x-y,y)$ as on left of Fig.3. Substituting this into Eq.10 the Eq.2 is recovered as $c(x+y,y) = (x+y)(x+y+1)/2 + y = C_2(x,y)$. The following definition (term) will be used below for this kind of algebraic derivation for vertical segments:

**Definition 3.** We call "algebraic derivation for vertical segments" to obtain the parameters for Cantor-like polynomials (like the case of right of Fig.2) if the mapping is "tiled" with vertical segments.

Triangle domains on Figs.2e-f have right angles like Fig.1, but simple $45°$ rotation from Figs.1 with Eqs.1-2 to Figs.2e-f to obtain the pairing polynomial is not possible, since segments on Figs.2e-f contain $1,3,5,7,\ldots 2n-1\ldots$ lattice points while diagonals on Fig.1 contain $1,2,3,4,5\ldots n\ldots$ ones. The bijective mapping on left of Figure 4 is the sum of math series (with $a_1=2$ and $d=2$) on x axis because $P_{\text{left of Fig.4}}(x= 0,1,2,3,4,5,y=0)= 0, 2, 2+4=6, 2+4+6=12, 2+4+6+8=20, 2+4+6+8+10=30$, etc., generally $P_{\text{left of Fig.4}}(x,y=0)= x(x+1)$, and at constant x the vertical y increases 1 by 1, so

$$P_{\text{left of Fig.4}}(x,y)= x(x+1)+y \quad \text{for } -x \leq y \leq x \quad \text{and} \quad P_{\text{right of Fig.4}}(x,y)= y(y+1)-x \quad \text{for } -y \leq x \leq y. \quad (13)$$

The second equation comes from analogous derivation or a simple $-90°$ rotation from left of Fig.4 to right of Fig.4, the latter is $(x,y) \rightarrow (y,-x)$. Eq.13 can also be derived by algorithm in Eq.7 with choosing e.g. the six points with function values $\{0,1,2,3,5,7\}$.

The out-of-domain behavior of Cantor-like polynomials is as follows: Only one example is exhibited which shows the general property in this respect. right of Fig.4 shows the out-of-domain values of Eq.13 along with the domain value. As defined, the Cantor domain of $P_{\text{right of Fig.4}}(x,y)$ is the lattice points at $y \geq 0$ and between the positive graphs of $y=\pm x$, the values without brackets. The values in brackets are the function values of $P_{\text{right of Fig.4}}(x,y)$ on lattice points positioned out-of-domain. All those are the scale (or values of) x (or y) axis or set Z shifted in relation to each other now horizontally, but generally, horizontally (e.g. right of Fig.4), vertically (e.g. right of Fig.2 and left of Fig.4) or diagonally (e.g. Fig.1 and left of Fig.2), that is:

**Theorem 1.** The tile property of Cantor-like polynomials is as follows: The graphs of Cantor-like polynomials can be considered as follows: Tile the 2D Cartesian (x,y) plain with many parallel horizontal (or vertical or diagonal) Z (or $N_0$) 1D lines (bands) shifted to each other with a systematically increasing value (the increasing sum of the generating math series $(2a_1+(n-1)d)/2$ via the "order number of tile" n), and the semi-infinite domain for Cantor-like polynomial is defined with the help of certain borderlines (e.g. $y=\pm x$, y=const., etc.) one cuts out from the 2D Cartesian plain. See right of Fig.4 as an example.

If the increase/decrease is vertical, then no $y^2$ term (e.g. left of Eq.13), if horizontal, then no $x^2$ term (e.g. right of Eq.13), and if diagonal, then linear combination of $2^{nd}$ order terms (e.g. Eqs.1-2) come up in Cantor-like polynomials; saw mapping in Eq.9 does not belong in this category, because that is defined on semi-infinite domain only in 1D. These will be stated more generally at the end in Theorem 7.

## EXTENSION OF DOMAIN AND THE BIJECTIVE MAPPING IN $N_0^2 \to N_0$

The known Rosenberg-Strong polynomials (R(x,y)) [4] are based on generation on left of Figure 5. In relation to Fueter-Pólya theorem, it cannot be one second order polynomial (as Eqs.1-2 on Fig.1), but must be two connected second order polynomials. The two connected polynomials come from the fact that in its generation there are "corners" or "line breaks" inside the domain. We provide little more details to explain this and to discover its properties. It provides a good example to see that fitting procedure in Eq.7 is much simpler than the derivation using subsequent mathematical series. Notice the math series on x axis {0, 0+1, 0+1+3, 0+1+3+5,...}, i.e. R(x,0)= [1+(2(x-1)+1)]x/2= $x^2$, as well as the math series on y axis {0, 0+3, 0+3+5, 0+3+5+7,...}, i.e. R(0,y)= [3+(2y+1)]y/2= $y^2+2y$, see left of Fig.5. Notice that every column between line y=x and axis x on left of Fig.5 the increment is one, i.e. y as in Eq.15, as well as between axis y and line y=x the increment is as R(0,y), a math series, R(x,y)= $x^2+x+(a_1+a_n)n/2$, where the series over x is {x+1,x+2,...,y}, i.e. $a_1=2x+3$ and $a_n=2(x+(y-x))+3$, i.e. $x^2+x + (x+y+2)(y-x) = y^2+2y-x$ as in Eq.14:

**Theorem 2.** The Rosenberg-Strong polynomials in relation to their joining lines follow. The bijective mapping R(x,y): $N_0^2 \to N_0$ defined on left of Fig.5 is

$$R(x,y) = y^2+2y-x \quad \text{if } 0\leq x<y \quad \text{(region I, left of Fig.5)} \qquad (14)$$
$$R(x,y) = x^2+y \quad \text{if } 0\leq y\leq x \quad \text{(region II, left of Fig.5)}. \qquad (15)$$

On line y=x the R(x,y=x)= $x^2+x= y^2+y$, that is, the two polynomials connect with different slopes. On left of Fig.5 the analogue/counter definition of this bijective mapping is the mirror image onto line y=x, Eqs.14-15 hold accordingly (via exchange x↔y).

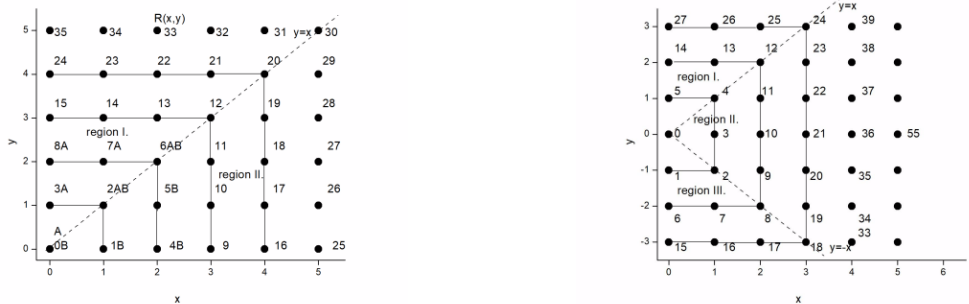

**FIGURE 5.** Bijective mapping of Rosenberg-Strong polynomials R(x,y): $N_0^2 \to N_0$ in Eqs.14-15, region I is between axis y and line y=x described by Eq.14 and region II is between line y=x and axis x described by Eq.15, both (including the border lines) are in quarter (+,+). For Eqs.14-15, two sets of six points with small coordinates have been chosen to find the parameters of polynomial in Eq.6 by the algorithm in Eq.7 marked with letters A and B in the two regions. (Derivation of Eqs.14-15 with using math series can be found in literature.) [Left]. Bijective mapping (half-square-spiral) of $P_{right\ of\ Fig.5}(x,y)$: {(x,y)| x∈$N_0$, y∈Z}⊂$Z^2 \to N_0$ in Eqs.16-18, regions I-II-III are separated by axis y and lines y=±x in half-plain (+,±) [Right].

Compact formulae using max(x,y) as counterpart of Eqs.14-15 can be found in ref.[4], however, Eqs.14-15 reflect better the broken nature of polynomials at line y=x. Another proof for Eqs.14-15 comes from the algorithm in Eq.7. One polynomial in Eq.6 is not adequate to fit (all points) for regions I-II together on left of Fig.5 by the reason just discussed, but two polynomials are necessary and have to be fitted separately. The y=x line belongs to both regions, for Eq.6 the six-six points with smallest values as possible should be picked for convenient calculation. (The chosen set of points should not be vertically or horizontally collinear to avoid zero determinant in the fit, as well as the set for fit must represent the periodicity). For example on left of Fig.5, set $\{(x,y,R)| (0,0,0), (0,1,3), (0,2,8), (1,1,2), (1,2,7), (2,2,6)\}$ marked with letter A in region I yields Eq.14, and set $\{(x,y,R)| (0,0,0), (1,0,1), (1,1,2), (2,0,4), (2,1,5), (2,2,6)\}$ marked with letter B in region II yields Eq.15 by using the algorithm in Eq.7. This algorithm is more systematic for the two regions and simpler procedure than the derivation with math series, and the same is true for more complex patterns below.

For the mappings on Figs.3b-5, the derivation of the polynomials is very hectic with math series, although elementary. However, with the algorithm in Eq.7 it is short, quick and convenient. The reason that algorithm in Eq.7 works for connected polynomials, like for Eqs.14-15 is that the "math series nature" preserves between regions on left of Fig.5, that is, the systematic increment is the corresponding additive constant among image values in the generation. To obtain the polynomials for bijective mapping on right of Fig.5 with algorithm in Eq.7, the coordinates with smallest image values can also be used: The sets of six points $\{(x,y,P)|(0,0,0), (0,1,5), (1,1,4), (0,2,14), (1,2,13), (2,2,12)\}$ in region I, $\{(x,y,P)|(0,0,0), (1,-1,2), (1,0,3), (1,1,4), (2,-2,8), (2,2,12)\}$ in region II and $\{(x,y,P)|(0,0,0), (0,-1,1), (1,-1,2), (0,-2,6), (1,-2,7), (2,-2,8)\}$ in region III yield the following polynomials.

**Theorem 3**: Bijective mapping of $P_{\text{right of Fig.5}}(x,y)$: $\{(x,y)|\ x\in N_0,\ y\in Z\}\subset Z^2\to N_0$ defined on right of Fig.5 is

$$P_{\text{right of Fig.5}}(x,y) = 2y^2-x+3y \quad \text{region I} \tag{16}$$
$$P_{\text{right of Fig.5}}(x,y) = 2x^2+x+y \quad \text{region II} \tag{17}$$
$$P_{\text{right of Fig.5}}(x,y) = 2y^2+x+y \quad \text{region III} \tag{18}$$

On border line $y=x$ the two polynomials in regions I-II connect with different slopes, so the two polynomials in regions II-III on border line $y=-x$. On right of Fig.5 the analogue definition of this bijective mapping is the mirror image onto axis x, Eqs.16-18 hold accordingly (via $y\leftrightarrow -y$).

The proof follows. Alternatively to the use of algorithm in Eq.7, algebraic derivation with math series for vertical segments to obtain Eq.17 for region II can be found in App.2, etc..

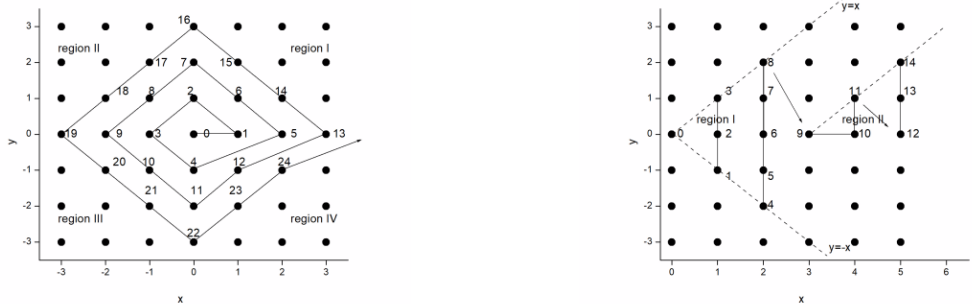

**FIGURE 6.** Bijective mapping (rhombus-spiral) of $P_{\text{left of Fig.6}}(x,y)$: $\{(x,y)|\ x\in Z,\ y\in Z\}=Z^2\to N_0$ in Eqs.19-22. Region I, II, III and IV are the (+,+), (-,+), (-,-) and (+,-) Cartesian quarters resp. [Left]. Bijective mapping (a complex domain) of $P_{\text{right of Fig.6}}(x,y)$: $\{(x,y)|\ x\in N_0,\ y\in Z\}\subset Z^2\to N_0$ in Eqs.24-25. Region I ($0\leq|y|\leq x\leq 2$) and region II ($3\leq x<\infty$, $0\leq y\leq x-3$) are connected between points at (2,2) and (3,0) [Right].

Continuing the application of this simple algorithm in Eq.7, some more complex bijective mappings follow. The original and Cantor-like polynomials on Figs.1-3 have the common features that on the segments (or sides) mathematical series (with increment 1) take up space. It is quite obvious e.g. on right of Fig.2 at a chosen $x_0$ that, Eq.10 gives $c(x_0,y)= c(x_0,0)+y$ in which the y changes one by one. This is the reason that $2^{nd}$ order polynomials cannot describe breaks like the ones on Fig.3, and connected polynomials are needed. Bijective mapping on left of Fig.6 shows the important feature of the borderlines of domains discussed with Eq.7, here between 4→5, 12→13, etc., and indeed the sides in region IV are not math series ($\{11,12,5\}$, $\{22,23,24,13\}$, etc.. However, knowing that here the x axis does not belong to region IV the residues are math series as $\{11,12\}$, $\{22,23,24\}$, etc.. In this way the regions are joined on left of Fig.6 as follows: Region I-II and II-III and III-IV have common lattice points on +y and -x and -y axis, respectively. However, region IV-I do not have common lattice points on +x axis, all lattice points on x axis belong to region I. There is another interesting feature, the point in origin $P(0,0)=0$ does not belong to any region. The reason is that its slope (or jump) to the next rhombus is different from others, compare jumps between rhombuses from image value jump P=0→P=1 (tg α=0) to jumps 4→5 and 12→13 and 24→25, etc. (tg α=½). These must be taken into account when polynomial in Eq.6 is fitted by the algorithm in Eq.7. For example, for region I the points with image values {6,14,15,26,27,28}, for region II the {8,17,18,30,31,32}, for region III the {10,20,21,34,35,36}, as well as for region IV the {12,23,24,38,39,40,} can be picked up for polynomial fit for Eq.6 by algorithm in Eq.7, but origin must not be included (otherwise false polynomials obtained). These points yield the following theorem.

**Theorem 4**: Bijective mapping of $P_{\text{left of Fig.6}}(x,y)$: $\{(x,y)| x \in Z, y \in Z\} = Z^2 \rightarrow N_0$ defined on left of Fig.6 is $P_{\text{left of Fig.6}}(0,0) = 0$,

$P_{\text{left of Fig.6}}(x,y) = 2x^2 + 4xy + 2y^2 - 2x - y + 1$      region I including +x and +y axis except origin      (19)
$P_{\text{left of Fig.6}}(x,y) = 2x^2 - 4xy + 2y^2 \quad -y + 1$      region II including +y and –x axis except origin      (20)
$P_{\text{left of Fig.6}}(x,y) = 2x^2 + 4xy + 2y^2 \quad -y + 1$      region III including –x and –y axis except origin      (21)
$P_{\text{left of Fig.6}}(x,y) = 2x^2 - 4xy + 2y^2 + 2x - y + 1$      region IV including –y axis but not +x axis.      (22)

Notice that in the four regions the Fig.1 with Eqs.1-2 is transformed by Lemmas 2-3. As a proof, use algorithm in Eq.7, e.g. pick the six points with coordinates (near to origin) for the polynomial fit which represents the periodicity correctly as mentioned before the theorem. Alternative proof with math series include reasoning as the concentric rhombuses have 4n-4 lattice points, so on x axis the P(x,y=0) polynomial has this increment (region I), i.e. P(x>1,0)=P(x-1,0)+4x-4, etc..

The two corresponding polynomials in regions I-II, II-III and III-IV connect in cyclic way with different slopes. On left of Fig.6 the natural numbers run increasingly in counter-clockwise direction, the clockwise direction, or the mirror image to x or y axis can be calculated in the same way as for Eqs.19-22 or can be transformed as in case of Theorem 3. (If one fits Eq.6 to six points from region IV and axis x, e.g. with image values {0,4,1,11,12,5} on left of Fig.6, the polynomial works for these six points, but gives false values for other lattice points in region IV, etc..) Eqs.19-22 can be written in one equation with logical functions as

$P_{\text{left of Fig.6}}(x,y) = 2x^2 + 4\operatorname{sgn}(x)\operatorname{sgn}(y)xy + 2y^2 - 2H(x)\operatorname{sgn}(y)x - y + 1$      except origin.      (23)

Connecting different patterns of domains is also possible like on right of Fig.6, wherein finite beginning part on left of Fig.4 and semi-infinite region on right of Fig.2 are connected:

**Theorem 5.** Bijective mapping of $P_{\text{right of Fig.6}}(x,y)$: $\{(x,y)| x \in N_0, y \in Z\} \subset Z^2 \rightarrow N_0$ defined on right of Fig.6 is

$P_{\text{right of Fig.6}}(x,y) = P_{\text{left of Fig.4}}(x,y) = x(x+1) + y$      region I ($0 \leq |y| \leq x \leq 2$)      (24)
$P_{\text{right of Fig.6}}(x,y) = \frac{1}{2}(x-3)^2 + \frac{1}{2}(x-3) + y + 9$      region II ($3 \leq x < \infty$, $0 \leq y \leq x-3$).      (25)

As a proof, use algorithm in Eq.7 for Eq.24 with picking adequate six points of image values for the fit or recognize that it is the same as $P_{\text{left of Fig.4}}$ in Eq.13 in the beginning part. Use Lemma 2 and Eq.10 to obtain $P_{\text{right of Fig.6}}(x,y) = c(x-3,y)+9$ in Eq.25.

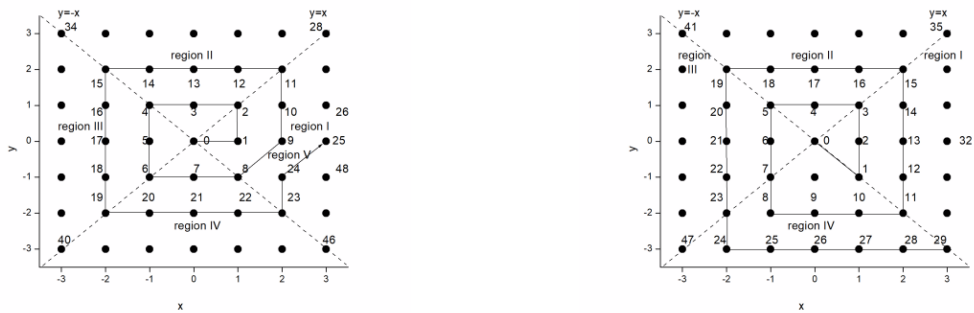

**FIGURE 7.** Bijective mapping (square-spiral on [Left] and rectangle-spiral on [Right]) of $P_{\text{left of Fig.7}}(x,y)$ and $P_{\text{right of Fig.7}}(x,y)$: $\{(x,y)| x \in Z, y \in Z\} = Z^2 \rightarrow N_0$ in Eqs.26-34

Topologically interesting how domains can be tiled, for example on Fig.5 called square-spiral (left of Fig.7) and rectangle-spiral (right of Fig.7). The domains for algorithm in Eq.7 must be defined such as the horizontal or vertical sides must contain finite math series with increment 1. For example, in region IV on left of Fig.7 the {6,7,8} or {19,20,21,22,23} are adequate, so lattice points on y=±x (and the internal ones) belong to the polynomial describing region IV. However, in region IV on right of Fig.7 the {23,8,9,10,11} is not a finite math series segment, only the {8,9,10,11}, so lattice points on y=-x (and internal ones) belong to region IV, but no lattice points on y=x belonging to region IV. The latter is indicated clearly by the connecting lines of image values between n and n+1. Notice that point (0,-1) on right of Fig.7 does not have image value, so more precisely, domain of $P_{\text{right of Fig.7}}$ is not $Z^2$ but $Z^2 \setminus \{(0,-$ 1)}. On Fig.5 the natural numbers from set $N_0$ run increasingly in counter-clockwise direction, the clockwise direction, or the mirror image to x or y axis can be calculated analogously. On lines $y=\pm x$ the two corresponding neighboring polynomials connect with same function value but different slopes if the borderline is common, e.g. regions I-II on left of Fig.7; (if no common points as between regions III-IV on right of Fig.7, this relation is irrelevant). Theorem 6 holds for Fig.5 with proof as using algorithm in Eq.7 and for simplicity picking up six points in the region in question near to the origin:

**Theorem 6.** Bijective mapping of $P_{\text{left of Fig.7}}(x,y)$: $\{(x,y)|\ x\in Z, y\in Z\}=Z^2\to N_0$ defined on left of Fig.7 is $P_{\text{left of Fig.7}}(0,0)=0$ and

$$P_{\text{left of Fig.7}}(x,y) = 4x^2\quad -4x + y +1 \qquad \text{region I}\ (0\leq y \leq x,\ x>0) \qquad (26)$$
$$P_{\text{left of Fig.7}}(x,y) = \quad 4y^2 -x -2y +1 \qquad \text{region II}\ (0\leq |x|\leq y,\ y>0) \qquad (27)$$
$$P_{\text{left of Fig.7}}(x,y) = 4x^2\quad\quad -y +1 \qquad \text{region III}\ (0\leq |y|\leq |x|,\ x<0) \qquad (28)$$
$$P_{\text{left of Fig.7}}(x,y) = \quad 4y^2 +x -2y +1 \qquad \text{region IV}\ (0\leq |x|\leq |y|,\ y<0) \qquad (29)$$
$$P_{\text{left of Fig.7}}(x,y) = 4x^2\quad +4x +y +1 \qquad \text{region V}\ (0<|y|\leq x,\ y<0) \qquad (30)$$

Bijective mapping of $P_{\text{right of Fig.7}}(x,y)$: $\{(x,y)|\ x\in Z, y\in Z\}=Z^2\setminus\{(0,-1)\}\to N_0$ defined on right of Fig.7 is $P_{\text{right of Fig.7}}(0,0)=0$, the $P_{\text{right of Fig.7}}(0,-1)$ is not defined and

$$P_{\text{right of Fig.7}}(x,y) = 4x^2\quad -x +y -1 \qquad \text{region I}\ (\text{between and on } y=\pm x,\ \text{as well as } x>0) \qquad (31)$$
$$P_{\text{right of Fig.7}}(x,y) = \quad 4y^2 -x +y -1 \qquad \text{region II}\ (\text{between and on } y=\pm x,\ \text{as well as } y>0) \qquad (32)$$
$$P_{\text{right of Fig.7}}(x,y) = 4x^2\quad -3x -y -1 \qquad \text{region III}\ (\text{between and on } y=\pm x,\ \text{as well as } x<0) \qquad (33)$$
$$P_{\text{right of Fig.7}}(x,y) = \quad 4y^2 +x +3y -1 \qquad \text{region IV}\ (\text{between } y=\pm x\ \text{and on } y=-x,\ \text{as well as } y<-1) \qquad (34)$$

The proof follows. To obtain the polynomials with applying math series is hectic in this case, but elementary. Use the convenient algorithm in Eq.7 for Eqs.26-34, we list the adequate sets of six points (with low value coordinates) for representation. For Eqs.26-30 the choice of points is as follows: Image value 0 cannot be included for any region by the same reason as on left of Fig.6, i.e. the slope (jump from n to n+1) at $0\to 1$ is not the same as at $8\to 9$, $24\to 25$ and $48\to 49$, etc.. For example, for region I the points with image values $\{10,26,27,50,51,52\}$ yielding $\det[A_{6x6}]=-4$, for region II the $\{3,12,13,29,30,31\}$ and notice that pattern with $45^0$ is used for domain with $90^0$, as well as the $\det[A_{6x6}]=-4$, for region III the $\{5,17,37,16,36,35\}$, for region IV the $\{7,20,21,41,42,43\}$ and for region V the $\{24,47,48,78,79,80\}$ can be picked up for polynomial fit for Eq.6 by algorithm in Eq.7 yielding Eqs.26-30. For Eqs.31-34 the choice of points is as follows: Image value 0 cannot be included for any region by the same reason as on left of Fig.6, the angle (jump $n\to n+1\to n+2$) at $0\to 1\to 2$ is $45^0$ while if value n+1 with n>0 is on $y=\pm x$, the angle of jump is $90^0$. For example, for region I the points with image values $\{2,13,14,32,33,34\}$ yielding $\det[A_{6x6}]=-4$ or alternatively $\{2,12,13,14,31,33\}$ yielding $\det[A_{6x6}]=-8$, for region II the $\{4,16,17,36,37,38\}$, for region III the $\{6,20,21,42,43,44\}$ and for region IV the $\{9,26,27,51,52,53\}$ can be picked up for polynomial fit for Eq.6 by algorithm in Eq.7 yielding Eqs.31-34. Alternatively, algebraic derivation for vertical segments for Eq.26 for region I on left of Fig.7 can be found in App.3.

Regions I and II of polynomial $P_{\text{right of Fig.7}}(x,y)$ connect at $y=x$, where $P_{\text{right of Fig.7}}(x,y=x)=4x^2-1$, while between region III and IV they do not. Figs.4a and 5a-b are drawn in spiral way, but those can be considered as concentric rhombuses, squares and rectangles, resp., as drawn on Fig.3. Cantor-like polynomials for other mappings than exhibited on figures in this work can be evaluated accordingly.

Finally, the following two theorems can be stated:

**Theorem 7.** If the direction (graph) of the increasing integers from image set (generally $N_0$) in the domain (full or subset of $Z^2$) is on diagonal segments, the Cantor-like polynomials contain all second order terms, i.e. $x^2$, $xy$ and $y^2$ (as well as first order term(s)), for example, Figs.1, 2a and 4a with the corresponding Eqs.1-2, 5 and 19-23. If the direction is vertical (or horizontal), the Cantor-like polynomials contain second order term $x^2$ (or $y^2$) only (as well as first order term(s)), for example, Figs.2b, 2e-2f, 3 and 5 with the corresponding Eqs.10, 13, 14-18 and 22-30. In the previous statement the domain is semi-infinite (e.g. $N_0^2$ or the lattice points of (+,+) quarter of Cartesian plain on Fig.1) or infinite (e.g. the $Z^2$ or the lattice points of full Cartesian plain on Fig.5) in "two dimension" (i.e. contains more than one non-parallel semi-infinite or infinitely long straight line through lattice points (e.g. $y=nx$ with $n\in Z$ on Fig.5)); if the domain is semi-infinite in "one dimension" (i.e. contains only parallel semi-infinite long straight lines, e.g. $y=k$ with $k=0,1,\ldots n-1$ in Table 1) the Cantor-like polynomials are first order only.

**Theorem 8.** All Cantor-like polynomials have a corresponding Diophantine equation of which properties can be described in relation to the full set natural numbers. For example, see Lemma 1 above or the Diophantine equation to Eq.10 for right of Fig.2 is $2c = a^2 + a + 2b$, and this equation for $a=0,1,2,\ldots,\infty$, $0\leq b\leq a$ picks up 1., all natural numbers from $N_0$ and 2, once and only once.

The proof of example in Theorem 8 comes simply from the way we arrived to Eq.10 (so for other cases), as well as the reason that all $(a,b)\in N_0^2$ yield even number on the left providing integer c (i.e. integers on lattice points) is that if "a" is even or odd, the $a^2 + a$ is always even (because $(2A)^2 + (2A) = 2(2A^2 + A)$ is even or $(2A+1)^2 + (2A+1) = 2(2A^2+3A+1)$ is also even).

This theorem looks simple, but far not trivial, e.g. leave $a^2$ from the equation: $2c=a+2b$ yields half-integer $c=3/2$ only if $(a,b)=(3,0)$ or at $(a,b)=(2,1)$ and $(4,0)$ it yields $c=2$ for both, i.e. not bijective.

Further manipulations are also possible as follows. Mirror the $(+,\pm)$ half-plain on right of Fig.5 (Eqs.16-18) to axis y and shift it to left by 1 ($P_{\text{right of Fig.5}}(-1-x,y)$), make the image values on $(+,\pm)$ half-plain even as $2P$ and make the image values on $(-,\pm)$ half-plain odd as $2P+1$, the

$$H(-x)\,(2\,P_{\text{right of Fig.5}}(-1-x,y) + 1) + 2\,H(x+\varepsilon)\,P_{\text{right of Fig.5}}(x,y) : Z^2 \to N_0 \quad (35)$$

($\varepsilon>0$ any real number) bijectively maps the entire $Z^2$ to the entire $N_0$ in a "very" zig-zag way e.g. in comparison to Figs.1-5. On Fig.1 the distance between adjacent image values, $d(n,n+1)$, is $\sqrt{2}$ or the increasing length of diagonals between points on axis x and y at "corners" by Pythagorean theorem, similar holds for Eq.35. Instead of $2P$ and $2P+1$ used in Eq.35 the $P$ (unchanged) and $-P-1$ resp. makes bijective mapping $Z^2 \to Z$ using full $Z^2$ and full $Z$.

**TABLE 1.** Pairing (Cantor-like) polynomials $P(x,y)$ when $x \in N_0$ but $y \in [0,n-1]$ is bounded.

| Mapping | | | | | | | shape | y width | Pairing polynomial | Pairing polynomial if width is n |
|---|---|---|---|---|---|---|---|---|---|---|
| | 1 | 3 | 5 | 7 | 9 | 11 | saw | 2 | $P = 2x-y$ | $P = nx-(n-1)y$ |
| 0 | 2 | 4 | 6 | 8 | 10 | 12 | | | $x \in N_0,\ y \in [0,1]$ | $x \in N_0,\ y \in [0,n-1]$ |
| | | | | 4 | 9 | 14 | saw | 5 | $P = 5x-4y$ | $P = nx-(n-1)y$ |
| | | | 3 | 8 | 13 | 18 | | | $x \in N_0,\ y \in [0,4]$ | $x \in N_0,\ y \in [0,n-1]$ |
| | | 2 | 7 | 12 | 17 | 22 | | | | |
| | 1 | 6 | 11 | 16 | 21 | 26 | | | | |
| 0 | 5 | 10 | 15 | 20 | 25 | 30 | | | | |
| 2 | 5 | 8 | 11 | 14 | 17 | 20 | comb | 3 | $P = 3x+y$ | $P = nx+y$ |
| 1 | 4 | 7 | 10 | 13 | 16 | 19 | | | $x \in N_0,\ y \in [0,2]$ | $x \in N_0,\ y \in [0,n-1]$ |
| 0 | 3 | 6 | 9 | 12 | 15 | 18 | | | | |

The $\sin^2((x+y)\pi/2) = \{1\text{ or }0\}$ and $\cos^2((x+y)\pi/2) = \{0\text{ or }1\}$ if $x+y = \{\text{odd or even}\}$, resp.. With these the two cases on Fig.1 can be mixed as alternating the direction of image values $n \to n+1$ on adjacent lines $x+y = \{2N\text{ and }2N+1\}$ by the 2-fold polynomial using Eqs.1-2 as

$$P_{\text{left of Fig.8}}(x+y) = C_{1\text{or}2}(x,y)\sin^2((x+y)\pi/2) + C_{2\text{or}1}(x,y)\cos^2((x+y)\pi/2) \quad (36)$$

In relation to Fueter-Pólya theorem, substituting the Taylor power expansion of sinus and cosinus functions into Eq.36, it becomes other bijective polynomial than Eqs.1-2, but not a $2^{nd}$ order polynomial (but infinite).

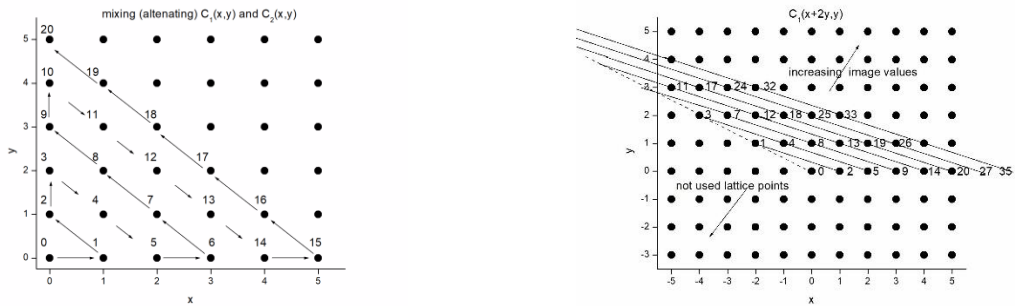

**FIGURE 8.** Bijective mappings $P_{\text{left of Fig.8}}(x,y): N_0^2 \to N_0$ [Left] and $P_{\text{right of Fig.8}}(x,y): \{(x,y)|\ x \in Z$ and $y \geq 0$ and $y \geq -x/k$ and $k=2\} \to N_0$ in Eq.36 ($C_1 \sin^2(.) + C_2 \cos^2(.)$) and Eq.37, resp. [Right].

Extension of the domain of original Cantor polynomials in Eqs.1-2 on Fig.1 is possible with the general trick on left of Fig.3, for example, shift every horizontal strip horizontally on left of Fig.1 to left by $k \in N_0$ steps, i.e.

$$P_{\text{right of Fig.8}}(x,y) = C_1(x+ky, y) = \tfrac{1}{2}((x+ky)^2 + 2(x+ky)y + y^2 + 3(x+ky) + y) \,. \tag{37}$$

The k=2 case is plotted on right of Figure 8. Alternatively, $C_2(x,y)$ can also be used. The domain is extended to the larger $\{(x,y)|\ x \in Z$ and $y \geq 0$ and $y \geq -x/k\}$ containing $N_0^2$. The important thing is that one polynomial is enough in Eq.37 on this larger domain, no need for connected polynomials as e.g. on Figs.3-5, but indeed no corners and fault lines. The domain tends to $Z \times N_0$, i.e. to $(\pm,+)$ half-plain, but the image values 1,3,6, etc. get further and farther from the origin.

## EXTENSION OF CANTOR'S BIJECTIVE MAPPING TO $N_0^3 \to N_0$

Here comes the "gem" of the article. Table 2 caption explains this bijective mapping. The $P_{3D}(x,y,z) = a_{300}x^3 + a_{210}x^2y + \ldots a_{003}z^3 + \ldots a_{100}x + a_{010}y + a_{001}z + a_{000}$ is 20 parameters 3rd order 3 variable polynomial. The number of points in "initial tetrahedron" for analogue algorithm in Eq.7 is 1+3+6+10=20 (listed in Table 1 plus the origin). The linear system can be solved block-wise (in its 20x20 matrix), for example, $P_{3D}(0,0,0)=0$ yields $a_{000}=0$ and parameters in $P_{3D}(0,0,z) = a_{003}z^3 + a_{002}z^2 + a_{001}z$ can be obtained via z=1,2,3 with $P_{3D}=1, 4, 10$, etc..

**TABLE 2.** Parallel triangles (cross sections) with corners (0,0,N), (N,0,0) and (0,N,0) perpendicular to line x=y=z in 3D Cartesian space in its (+,+,+) eight part. The 3 digit numbers notates the coordinates (x,y,z) of lattice points, e.g. 001=(0,0,1), and 1 digit numbers are the polynomial values $P_{3D}(x,y,z)$. Triangle-0 is degenerated to point origin for $P_{3D}(0,0,0)=0$. Image numbering starts on the (N,0,0) corners of the next triangle, and continued as decreasing the z coordinates "1 by 1" and from axis x to axis y at constant z, another orientation is from axis y to axis x, all together (cyclic change of x,y and z) 3!=6 equivalent polynomials exist for $P_{3D}(x,y,z)$. Notice that all triangles have math series as on Fig.1. This is a bijective mapping $P_{3D}(x,y,z): N_0^3 \to N_0$, the generalization of right of Fig.1 from 2D ($C_2(x,y): N_0^2 \to N_0$, Fig.1 containing 2!=2 equivalent polynomials in Eqs.1-2) to 3D,. The lines y=-x+N as x+y=N on right of Fig.1 correspond to Triangle-N on plain x+y+z=N here.

| Triangle-1 (x+y+z=1) 3 lattice points | | | | Triangle-2 (x+y+z=2 for coordinates) 6 lattice points | | | | | | Tringle-3 (x+y+z=3 for coordinates) 10 lattice points | | | | | |
|---|---|---|---|---|---|---|---|---|---|---|---|---|---|---|---|
| | 001 1 | | | | | 002 4 | | | | | | | 003 10 | | |
| 100 2 | | 010 3 | | | 101 5 | | 011 6 | | | | 102 11 | | 012 12 | | |
| | | | | 200 7 | | 110 8 | | 020 9 | | 201 13 | | 111 14 | | 021 15 | |
| | | | | | | | | | | 300 16 | | 210 17 | | 120 18 | 030 19 |

After bijective $P_{1D}(x)=x: N_0 \to N_0$, $P_{\text{squares}}(x)=x^2: N_0 \to \{\text{integer squares}>0\}$ and in Eqs.1-2 the $C_{1\text{or}2}: N_0^2 \to N_0$, finally:

**Theorem 9.** The bijective $N_0^3 \to N_0$ polynomial (cyclic change of x,y and z makes it 6-fold) is

$$P_{3D}(x,y,z) = [x^3+y^3+z^3 + 3(xz^2+yz^2+zx^2+2xyz+zy^2+yx^2+xy^2) + 3(2x^2+2y^2+z^2+2xz+2yz+4xy) + 5x+11y+2z]/6 \,. \tag{38}$$

The generalization of bijective mappings y=x (k=1), Fig.1 (k=2) and Table 2 (k=3) is a k!-fold Cantor-like polynomial: $N_0^k \to N_0$ which runs on adjacent hiper-plains $\Sigma_{i=1}^k x_i = N \in N_0$ and starts with $(1/k!)\Sigma_{i=1}^k x_i^k$.

## CONCLUSIONS

The extension of Cantor polynomials have been described and analyzed with exciting eye-catching examples.

## APPENDIX

**Appendix 1.** The Heaviside step function (or the unit step function) is $H(x \leq 0) = 0$ and $H(x>0) = 1$, int(z) is the integer part or floor of z, e.g., int(5/2)= 2, $N_0$ is the natural numbers with zero, R(x,y) is the Rosenberg-Strong polynomials, sqrt(z)= $\sqrt{z}$, the sign function is sgn(x>0)=1, sgn(x<0)=-1 and sgn(0)=0, as well as Z is positive and negative integer numbers with zero.

To the saw mapping in Eq.9 plotted on right of Fig.3 some more mappings of similar kind (with finite width or "one dimension" semi-infinite domain) are listed. To save space these will be listed in more compact forms in Table 1 and discussed here. Saw mapping with width 2 is basically the vertical separation of odd and even integers. The set of points with function values {1,2,3,4,5,6} is adequate for algorithm in Eq.7 if width is 2, while set of points for saw mapping with width 5 containing function values {0,1,2,5,6,10} is adequate for algorithm in Eq.7, for example. The (x,y,P)= (0,0,0) point can be in the fit in saw mapping with width n, except when n=2, in the latter it yields zero determinant, but the linear system is not contradictive. More, if n=2, a more general pairing polynomial is P= $(2+a)y^2+2x-(3+a)y$ holds for saw mapping which is a trivial case because only y=0 and 1 is allowed in the domain. Notice that, the domain of original Cantor polynomials in Eqs.1-2 plotted on Fig.1 is semi-infinite in two direction (x and y) and the polynomial is 2$^{nd}$ order, while Table 1 contains Cantor-like polynomials with semi-infinite domain in one direction (especially x) and the polynomial is 1$^{st}$ order (see Theorem 7 in more detail).

**Appendix 2.** Function or image values along the perimeter of half-squares divided to quarter-squares by axis x on right of Fig.5 contains 2x+1 lattice points in quarter (+,+) including the positive x and y axis, e.g. at (1,0) the image value is 3 having 3 points (the image values 3,4,5), at (3,0) the image value is 21 having 7 points (the image values 21,…24,…27). If image value P(x-1,0) is known, at (x,0) the (+,+) quarter contains a=(2(x-1)+1)-1 and the (+,-) quarter contains b=2x+1 lattice points (along quarter-squares) to reach P(x,0) value, that is, P(x,0)= P(x-1,0)+a+b= P(x-1,0)+4x-1. (Notice that, e.g. from P(x=2,y=0)=10 to P(3,0)=21 the way is along two quarters of subsequent square-perimeters: 10,11,…14,15,…20,21.) In this way, see right of Fig.5, P(1,0)=3, P(2,0)= P(1,0)+4*2-1= 10= 3+7, P(3,0)= P(2,0)+4*3-1= 21= 3+7+11, P(4,0)= P(3,0)+4*4-1= 36= 3+7+11+15, P(5,0)=55, etc.. Finally, P(x,0)= $a_1+a_2+…a_x$= $x(a_1+a_x)/2$= $x(2a_1+(x-1)d)/2$= $2x^2+x$ because $a_1$=3 and d=4. The P(x,y)= $2x^2+x+y$ simply, since on vertical lines the y increases 1 by 1 in region II, as in Eq.17. The systematic algorithm in Eq.7 yields the same, but that can be handled more conveniently than the derivation with math series right in this Appendix.

**Appendix 3.** Function or image values along perimeter of squares crossing axis x at (x,0) on left of Fig.7 are distributed at 4(2x+1)-4=8x lattice points, e.g. at (3,0) the image values is 25 having 8*3=24 points (with values 25,…,28,…,34,...,40,…46,47,48). If image value P(x-1,0) is known, P(x,0)= P(x-1,0) + 8(x-1), because lattice points on perimeter of square including point (x-1,0) contain 8x lattice points added. In this way, see left of Fig.7 region I, P(1,0)=1, P(2,0)= P(1,0)+8*1= 1+8, P(3,0)= P(2,0)+8*2= 1+8+16, P(4,0)= P(3,0)+8*3= 1+8+16+24, P(5,0)= 1+8+16+24+32= 81, etc.. Finally, P(x,0)= 1+8(1+2+3+4+…(x-1))= 1+4x(x-1). The P(x,y)= P(x,0)+y= 1+4x(x-1)+y simply, since on vertical lines the y increases 1 by 1 in region I (0≤y≤x) on left of Fig.7, as in Eq.26. Again, systematic algorithm in Eq.7 is much simpler for Eq.26.

## ACKNOWLEDGMENTS

Financial and emotional support for this research from OTKA-K-NINCS 2015-115733 and 2016-119358 is kindly acknowledged. Thanks to Szeger Hermin for her help in typing the manuscript. The subject has been presented in ICNAAM_2022, Greece, and Am. Inst. Phys. Conf. Proc. 2023.

## REFERENCES

1. G. Cantor: J. Reine Angew. Math., Band **84**, 242–258 (1878)
2. R. Fueter, G. Pólya: Vierteljschr. Naturforsch. Ges. Zürich **68**, 380–386 (1923)
3. M. B. Nathanson: The Am. Math. Monthly, **123** (10), 1001-1012 (2016)
4. A. L. Rosenberg, H. R. Strong: IBM Technical Disclosure Bulletin, **14** (10), 3026–3028 (1972)
5. S. Kristyan: Am. Inst. Phys. Conf. Proc. **1863**, 560013 (2017)
6. S. Kristyan: Am. Inst. Phys. Conf. Proc. **1978**, 470064 (2018)